\documentclass[12pt]{amsart}
\usepackage{amsmath,amsthm,latexsym,amsfonts,amssymb}
\usepackage{datetime, graphicx}

\numberwithin{equation}{section}

\title {Limit cycles of planar vector fields: Hilbert's 16th problem and o-minimality}

\author {Patrick Speissegger}

\address {Department of Mathematics and Statistics, McMaster University, 1280
	Main Street West, Hamilton, Ontario L8S 4K1, Canada}

\email {speisseg@math.mcmaster.ca}

\date{\today\ at \currenttime}

\thanks{Supported by NSERC of Canada grant RGPIN 261961 and the Zukunftskolleg of Universit\"at Konstanz.  This note will appear in the Oberwolfach ``Snapshots of Modern Mathematics'' series.  I thank Zeinab Galal and Tobias Kaiser for vetting earlier versions of this note.}

\begin{document}
	
\maketitle
\markboth{Patrick Speissegger}{Limit cycles and o-minimality}

\begin{abstract}
	I discuss some recent work linking certain aspects of the second part of Hilbert's 16th problem to the theory of \hbox{o-minimality}.
	These notes are adapted from a lecture I gave in the Jour fixe seminar series at the Zukunfts\-kolleg of Universit\"at Konstanz in June 2017.
\end{abstract}

\section*{}

A \textit{vector field} is a map, denoted by $F$ below, that assigns to every point in the plane (or more generally in $n$-space) a \textit{vector}, which codes a \textit{direction} and a \textit{length}.  Two common examples of vector fields are \textit{force fields} (where the length represents acceleration) or \textit{fluid flow} (where the length represents speed).  Vector fields may or may not be time dependent; we assume here that $F$ is time independent.

The path traversed by an object following the vector field is called its \textit{trajectory}.  It is obtained by solving the differential equation $P'(t) = F(p(t))$, and it depends on the position $P(t_0)$ of the object at a given initial time $t_0$, called the \textit{initial condition}.  The branch of mathematics tasked with studying phenomena arising from vector fields is called \textit{dynamical systems} (see Perko \cite{Perko:2001fv} for an introduction).  Here is the fundamental theorem for trajectories of vector fields:\medskip

\noindent\textbf{Theorem \cite[Section 2.2]{Perko:2001fv}.}
\textsl{Under very mild assumptions on the vector field, each trajectory is uniquely determined by its initial condition.}
\medskip

Knowing the position $P(t_0)$ at time $t_0$ of an object following the vector field, and given a later time $t_1$, we would like to be able to predict the position $P(t_1)$ of this object at time $t_1$ (\textit{quantitative phenomena}) or give a general description of its long-term behaviour (\textit{qualitative phenomena}).   

For instance, if the vector field is linear, that is, $F(a) = Aa + b$, where $b$ is a fixed point in $n$-space and $A$ is an $n \times n$ matrix, then one can find explicit solutions in terms of the elementary functions $+$, $\cdot$, $\exp$, $\log$, $\sin$, $\cos$ and a few related functions and using complex numbers.  

However, for almost all other vector fields, no such explicit solutions exist.  One might say that the vocabulary of elementary mathematics is too small to describe the phenomena coded by vector fields.  In this sense, the goal of the field of dynamical systems can be stated as ``developing a mathematical vocabulary'' to describe such phenomena.
	
Quantitative phenomena are usually studied using Numerical Analysis, which is concerned with computing approximations to such phenomena.  This is done, for example, in models of weather or climate systems, of the burning processes in combustion engines, and many others.  However, numerical methods are not able to discern different qualitative phenomena, because the longer a model is run, the larger its inherent errors become. 
	
\subsection*{Some vocabulary to describe qualitative phenomena}

A \textit{singularity} of $F$ is a point $p$ such that $F(p) = 0$.  An object located at a singularity stays there forever; so its trajectory is a point.  The trajectory of an object is a \textit{cycle} (or \textit{periodic trajectory}), if the object revisits the same points periodically.  It follows from the theorem above that the trajectory of an object either visits every point at most once, or it is a cycle.  A trajectory \textit{spirals} if it is not a cycle but turns infinitely often around a fixed set of points called a \textit{limit set}.  

\begin{figure}[htbp]
	\begin{center}	
		\includegraphics[scale=0.2]{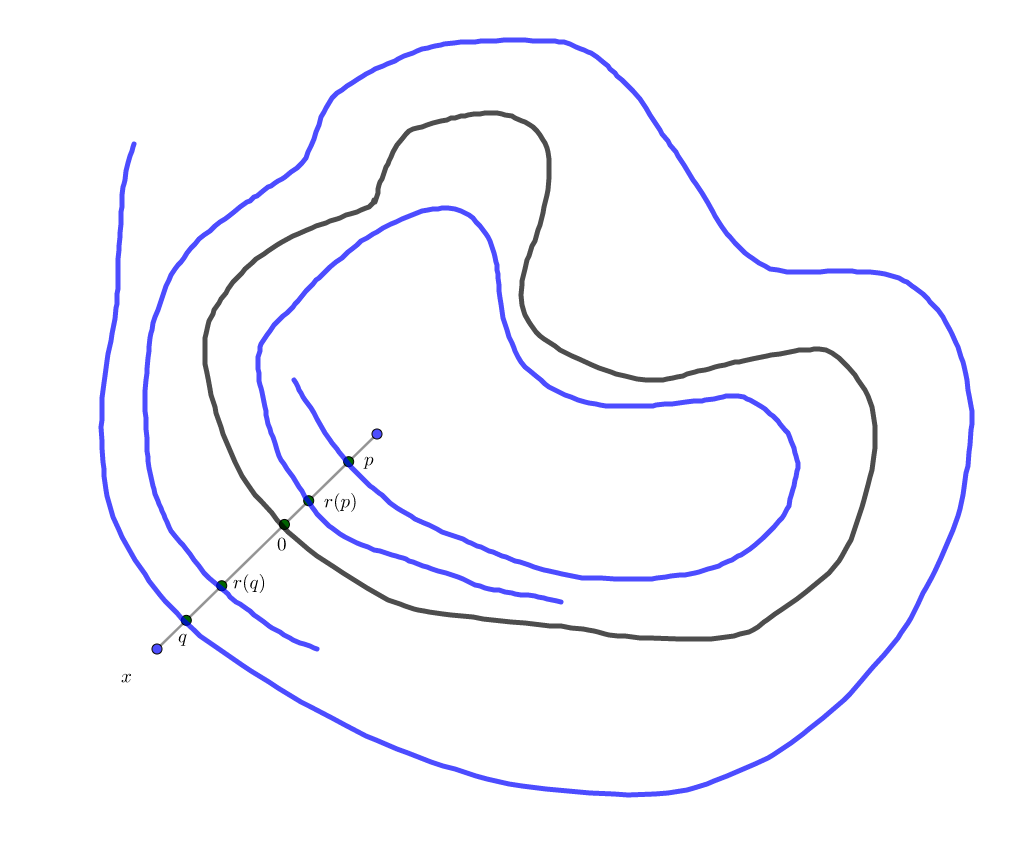}
	\end{center}
	\caption{\footnotesize A limit cycle (black) with nearby spiraling trajectories (blue).  The corresponding Poincar\'e first return map $r(x)$ is defined using a transverse segment.}
	\label{cap:limit}
	\end{figure}

A \textit{limit cycle} (see Figure \ref{cap:limit}) is a cycle such that every nearby trajectory spirals around it (in effect, either towards or away from it).  Thus, a limit cycle is the limit set of all trajectories near it.

	Predicting the qualitative phenomena of a time-independent vector field means answering the following question: if a particle is dropped at a specific point in $n$-space, then what is the nature of its trajectory?  To answer this question, we need to know where the singular points and the limit sets of the vector field are and, in the latter case, what their nature is (are all limit sets limit cycles, or are there other limit sets?).  A first step towards knowing the latter is to determine \textit{how many} singular points and limit cycles there are.
	
\subsection*{Counting singular points and limit cycles}

	We restrict our attention from now on to planar vector fields ($n=2$).\medskip

	\noindent\textbf{Example 1.} 		The planar vector field $F(x,y) = (y,-x)$ has one singularity, the origin.  All other trajectories are circles centered at the origin, hence cycles.  In particular, none of the trajectories are limit cycles.
	\medskip
	
	The vector field in the previous example is \textit{linear}, which is the simplest kind of vector field there is.  It is the only kind of vector field where we know how to count limit cycles--because there aren't any.
	
	Linear vector fields are examples of polynomial vector fields: 	a polynomial of degree $d$ in the variables $x$ and $y$ is an expression of the form $$a_{0,0} + a_{1,0}x + a_{0,1} y + a_{2,0} x^2 + a_{1,1} xy + \cdots + a_{d,0} x^d + \cdots + a_{0,d} y^d,$$  where the $a_{i,j}$ are real numbers. 
	\medskip

	\noindent Note: ``polynomial of degree 1'' is the same as ``linear''.  
	\medskip

	The vector field $F$ is polynomial of degree $d$, if each of the two components of $F$ is given by a polynomial of degree $d$. \medskip 
	
	\noindent\textbf{Example 2.}
		$F(x,y) = \left(x^2 + y^2, x-y^3\right)$ is polynomial of degree 3.
	\medskip
	
	After counting singularities and limit cycles of linear vector fields, which is easy, we could try counting singularities and  limit cycles of polynomial vector fields.  Since, in this case, singularities are just zeroes of polynomials, their study falls under the well-developed subject of real algebraic geometry, and I will not pursue it further here.  As to counting their limit cycles, the following was suggested by David Hilbert in his famous address given at the first  International Congress of Mathematics in the year 1900:\medskip
		
	\noindent\textbf{Hilbert's 16th problem (second part).}
	\textsl{If the vector field $F$ on the plane is polynomial of degree $d$, there exists a number $H(d)$ such that $F$ has at most $H(d)$ limit cycles.}
	
\subsection*{A very brief history of Hilbert's 16th problem}

	This problem remains open to this day.
	The timeline below is taken from Ilyashenko's more technical account \cite{Ilyashenko:2002fk} of Hilbert's 16th problem; I refer the reader there for  detailed references.\medskip
		
	\noindent\textbf{1923:}
	Dulac proves that every polynomial vector field has only finitely many limit cycles (\textbf{Dulac's problem}). 	His proof does not clarify if $H(d)$ exists, but his method proved to be useful for the study of dynamical systems in general. 
	\medskip
		
	\noindent\textbf{1955--57:}
			Petrovskii and Landis publish a solution of Hilbert's 16th problem.    It implies, in particular, that $H(2) = 3$. 
	\medskip
		
	\noindent\textbf{1963:}
			Ilyashenko and Novikov produce the first counterexamples to Petrovskii and Landis's solution (so their proof was wrong).
	\medskip
		
	\noindent\textbf{1979--80:}
			Chen, Wang and Shi give examples of quadratic (i.e., $d=2$) vector fields with 4 limit cycles; in particular, $H(2) \ge 4$.
	\medskip
		
	\noindent\textbf{1981:}
			Ilyashenko, in lectures given on Dulac's problem, discovers a previously overlooked gap in Dulac's proof.    In Ilyashenko's own words: ``Thus, after eighty years of development, our knowledge of Hilbert's 16th problem was almost the same as at the time when the problem was stated.'' 
	\medskip
		
	\noindent\textbf{1991--92:}
			Ecalle and Ilyashenko independently publish papers that fill the gap in Dulac's proof.    Both of these gap-filling proofs are much longer than Dulac's original proof, but show that Dulac's original argument was right in principle, ``just'' incomplete.
	\medskip
	
\subsection*{Poincar\'e's idea of how to count limit cycles}
	
	The idea is to reduce the two-dimensional counting problem (counting limit cycles in the plane) to a one-dimensional counting problem (counting certain points on a line). \medskip
	
	\noindent\textbf{Example 3 \cite[Section 3.4]{Perko:2001fv}:} to count limit cycles near a cycle (drawn in black in Figure \ref{cap:limit}), draw a line segment crossing the cycle that is not tangent to any trajectory.  Introduce a coordinate $x$ on this segment such that the intersection of the segment with the cycle is a $x=0$; so the line segment corresponds to an interval $(a,b)$ for some $a < 0 < b$.

	On this segment, define a map $r:(a,b) \longrightarrow (a,b)$ such that, for $x \in (a,b)$, the point $r(x)$ is the first intersection point of the trajectory going through $x$ with the segment that lies no farther to 0 than $x$; this map is called the \textit{Poincar\'e first return map}, see Figure \ref{cap:limit}.
	
	The point of this map is that counting cycles near a cycle corresponds to counting \textit{fixed points} of the associated map $r(x)$, i.e., points $x$ such that $r(x) = x$, near $x=0$.  Therefore, counting limit cycles near a cycle corresponds to counting \textit{isolated} fixed points of the map $r(x)$ near $x=0$, i.e., fixed points $x_0$ of $r(x)$ for which there exists an open interval about $x_0$ that contains no other fixed points of $r(x)$.
	
	The problem is that,
		while this reduces the dimension of the counting problem, it also takes us out of the realm of differential equations:  the Poincar\'e map $r(x)$ is not itself solution of any reasonably simple differential equation.  
		
	Poincar\'e overcame this problem by showing that the map $r(x)$ is \textit{analytic} at $x=0$, that is, it has a Taylor series expansion $$\hat{r}(x) = \sum_{n=0}^\infty a_n x^n = a_0 + a_1 x + a_2 x^2 + \cdots$$ at $x=0$, and this Taylor series \textit{converges}.  The latter implies
		that $r(x)$ can be approximately computed, to arbitrary precision, by computing a finite sum $a_0 + a_1 x + \cdots + a_n x^n$ for sufficiently large $n$.
	(For instance, all elementary functions mentioned earlier are analytic.)  The key observation about functions that are analytic at 0 is that their isolated fixed points cannot accumulate at 0.  Therefore:
	\medskip
	
	\noindent\textbf{Poincar\'e's corollary.}
	\textsl{The map $r(x)$ has only finitely many isolated fixed points near 0, so there are only finitely many limit cycles near a given cycle.}

\subsection*{Dulac's strategy for counting limit cycles}

\begin{figure}[htbp]
		\begin{center}
			\includegraphics[scale=.5]{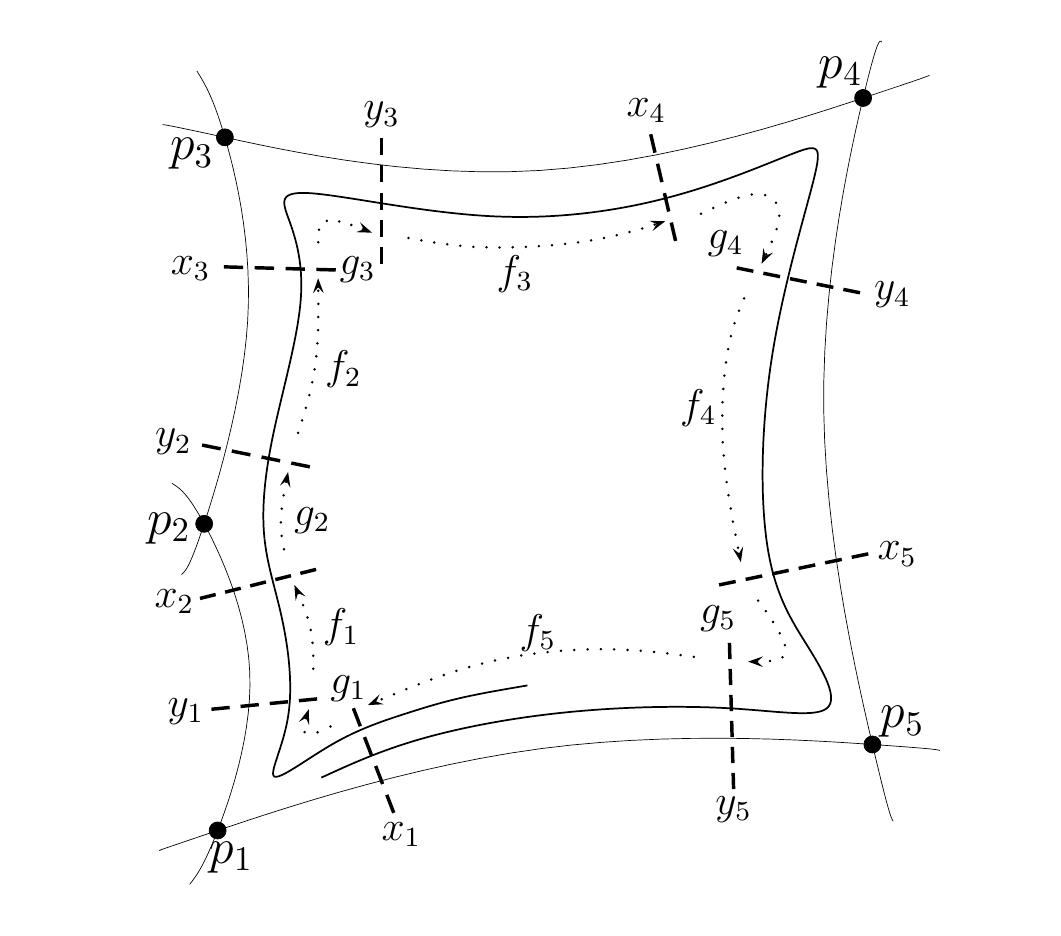}
		\end{center}
		\caption{\footnotesize A polycycle with associated transition maps $f_i$ and $g_i$}
		\label{cap:poly}
	\end{figure}

	Dulac showed that the general problem of counting limit cycles of polynomial vector fields can be reduced to a situation similar to that studied by Poincar\'e.  Here the cycle in Poincar\'e's situation is replaced by what is called a \textit{polycycle}, which is a closed curve consisting of finitely many singular points connected by trajectories as in Figure \ref{cap:poly}.  Using the transverse segment with coordinate $x = x_1$ in this figure, one can again define a corresponding first return map $r(x)$.  (The reason for the multiple segments in the figure will be explained later.)
	
	Thus, Dulac needed to prove that such a first return map $r(x)$ of a polycycle has finitely many isolated fixed points.  Similar to Poincar\'e's example, this can be done by showing the following:
	
	\begin{enumerate}
		\item these return maps $r(x)$ have \textit{asymptotic expansions} $\hat r(x)$ at $x=0$ (albeit more general than convergent Taylor series expansions);
		\item each such return map $r(x)$ is uniquely determined by its asymptotic expansion $\hat r(x)$.
	\end{enumerate} 
	
	While Dulac completed Point 1, Point 2 was the gap left unproved by him and proved 70 years later by Ecalle and Ilyashenko.
	
\subsection*{What else is needed for Hilbert's 16th problem?}

	For each degree $d$, let $\mathcal S_d$ be the collection of all polynomial vector fields in the plane of degree $d$.  Each vector field in $\mathcal S_d$ is denoted by $F_\mu$, where $\mu$ is a tuple of real numbers representing all the coefficients of the polynomials used in the definition of $F_\mu$. 
	
	To prove Hilbert's 16th problem, it is not enough to count the number of limit cycles near a polycycle of each vector field $F_\mu$ separately.  Instead, given a parameter $\mu$ and a polycycle $\Gamma$ of $F_\mu$, one needs to count all limit cycles near $\Gamma$ of all vector fields $F_{\mu'}$ for $\mu'$ close to $\mu$ (where ``close'' is to be understood in the sense of the usual topology on Euclidean spaces).  Moreover, more general \textit{limit periodic sets} (not defined here; it is sufficient for the purpose of this exposition to continue thinking of them as polycycles) need to be considered instead of polycycles.  Indeed, Roussarie \cite[Prop. 1 of Chapter 2]{MR3134495} shows that Hilbert's 16th problem follows if the following holds for every parameter $\mu$ and every limit periodic set $\Gamma$ of $F_\mu$:\medskip
	
	\noindent\textbf{Finite cyclicity conjecture (Roussarie).} \textsl{There exist a natural number $N$ and open neighborhoods $U$ of $\mu$ and $V$ of $\,\Gamma$ such that for every $\mu' \in U$, the vector field $F_{\mu'}$ has at most $N$ limit cycles contained in $V$.}\medskip
	
	Given a parameter $\mu$ and a limit periodic set $\Gamma$ of $F_\mu$, what makes the finite cyclicity conjecture difficult to prove (apart from the somewhat obscure nature of limit periodic sets in general) is that the return map $r_{\mu'}(x)$ of $F_{\mu'}$ around the limit periodic set $\Gamma$ is not necessarily well defined for all parameters $\mu'$ close to $\mu$ (because of so-called bifurcation phenomena, see \cite[Chapter 4]{Perko:2001fv}).  
	
	Assuming $\Gamma$ is a polycycle of $F_\mu$, one possible way to deal with this problem is to decompose $r_\mu(x)$ into the \textit{transition maps} $y_i = g_{\mu,i}(x_i)$ and $x_{i+1} = f_{\mu,i}(y_i)$ for $i=1, \dots, k$ as in Figure \ref{cap:poly}, where $k$ is the number of singularities on the polycycle $\Gamma$ (equal to 5 in the figure) and we convene that $x_{k+1} = x_1$.  One recovers the first return map from the transition maps as $$r_\mu(x) = (f_{\mu,k} \circ g_{\mu,k} \circ \cdots \circ f_{\mu,1} \circ g_{\mu,1})(x),$$ the successive composition of the $f_{\mu,i}$ and $g_{\mu,i}$.  
	By a general theorem on the dependence on initial conditions and parameters \cite[Section 2.3]{Perko:2001fv}, there are open neighbourhoods $U$ of $\mu$ and $V$ of $\Gamma$ such that the transition maps $f_{\mu',i}$ and $g_{\mu,i}$ are well defined for all parameters $\mu' \in U$ and segment coordinates $x_i, y_i \in V$ (although their composition may not be well defined if $\mu' \ne \mu$). 
	These parametric transition maps can be used, in place of the return maps, to describe the limit cycles of $F_{\mu'}$ near $\Gamma$: $x \in V$ corresponds to a limit cycle of $F_{\mu'}$ near $\Gamma$, with $\mu' \in U$, if and only if $x$ belongs to the set $A_{\mu'}$ of all isolated points of the set 
	\begin{multline*}
		\big\{x_1 \in V:\ \exists x_2, \dots, x_{k}, y_1, \dots, y_k \text{ such that } \\ y_i = g_{\mu',i}(x_{i}) \text{ and } x_{i+1} = f_{\mu',i}(y_i) \text{ for each } i\big\}.
	\end{multline*}   By the previous paragraph, the sets $A_{\mu'}$ are well defined for $\mu' \in U$, even if the composition of the transition maps is not well defined.  Similar parametric transition maps and corresponding sets $A_{\mu'}$ can be defined near every limit periodic set (not just polycycles).
	
	In Kaiser et al. \cite{Kaiser:2009ud}, this observation is used to formulate a criterion for these parametric transition maps that implies the corresponding finite cyclicity conjecture.  The new ingredient in this formulation comes from model theory, a branch of mathematical logic.

\subsection*{From Dulac's proof to Hilbert's 16th problem with \dots\ logic?}
	
	In the 1930s, G\"odel established some surprising implications of logic for the general study of mathematics, known as G\"odel's completeness and incompleteness theorems. 
	Out of these theorems, a new branch of mathematics called \textit{model theory} arose, started by Robinson in the 1950s.    
	It studies the implications of G\"odel's theorems for particular situations in mathematics; see Marker \cite{Marker:2002hq} for an introduction to model theory.   
	
	A crucial concept from model theory is that of \textit{definability}: given a set $\mathcal L$ of relations and functions on Euclidean space of various arities (the \textit{language}), we call \textit{$\mathcal L$-formula} any expression formed from symbols in $\mathcal L$, variables, the logical connectives $\wedge$ (``and''), $\vee$ (``or''), $\to$ (``implies'') and $\neg$ (``not''), as well as the logical quantifiers $\exists$ (``there exists'') and $\forall$ (``for all''), following the syntactic rules of first-order predicate logic \cite[Section 1.1]{Marker:2002hq}.  Of particular importance are the \textit{free} variables of an $\mathcal L$-formula $\phi$, that is, those variables in $\phi$ that are not bound by any quantifier in $\phi$.  A set $S \subseteq \mathbb R^n$ is \textit{definable} from $\mathcal L$ if there exists an $\mathcal L$-formula $\phi(x_1, \dots, x_n)$ with free variables among $x_1, \dots, x_n$ such that $$S = \left\{(a_1, \dots, a_n) \in \mathbb R^n:\ \phi(a_1, \dots, a_n) \text{ holds}\right\}.$$  The collection of all sets definable from $\mathcal L$ is referred to as an \textit{$\mathcal L$-structure (on the real numbers)}.  For instance, zerosets of polynomials are definable from the language $\mathcal L_{\text{or}} := \{+,-,\cdot, 0,1,=,<\}$ of ordered rings.  Of interest to this paper is the following:\medskip
	
	\noindent\textbf{Example 4.} If $\mathcal L$ is a language that contains $\mathcal L_{\text{or}}$ as well as the transition maps $f_i(\mu',x_i):= f_{\mu',i}(x_i)$ and $g_i(\mu',y_i):= g_{\mu',i}(y_i)$ defined on the set $U \times V$ for the polycycle $\Gamma$ above, then the set $$A:= \left\{(\mu',x) \in U \times V:\ x \in A_{\mu'}\right\}$$ is definable from $\mathcal L$.  This follows easily from the definition of $A$ above and the observation that $x$ is an isolated point of a set $S \subseteq \mathbb R$ if and only if there exists an $\epsilon > 0$ such that $S \cap (x-\epsilon,x+\epsilon) = \{x\}$.\medskip
	
	What makes this last example interesting in connection with Roussarie's conjecture is a tameness condition for $\mathcal L$-structures, now called \hbox{\textit{o-minimality}}, discovered by van den Dries and developed by Pillay and Steinhorn in the early 1980s; see van den Dries \cite{Dries:1998nj} for an introduction to o-minimality.  
	
	By definition, an $\mathcal L$-structure is \textit{o-minimal} if every subset of $\mathbb R$ definable from $\mathcal L$ is a finite union of intervals.  Since adding an existential quantifier to an $\mathcal L$-formula $\phi$ corresponds to taking a coordinate projection of the set defined by $\phi$, the collection of all sets definable from $\mathcal L$ is closed under taking coordinate projections.  
	Therefore, the o-minimality condition has implications for all sets definable from $\mathcal L$ (not just the subsets of $\mathbb R$).  In dimension greater than 1, the role of intervals is played by \textit{cells} in $\mathbb R^n$, which are defined by induction on $n$: the cells in $\mathbb R$ are the intervals and, if $n>1$, a cell in $\mathbb R^n$ definable from $\mathcal L$ is any of the sets $$\Gamma(f)_C:= \{(x,f(x):\ x \in C\},$$ $$(f,g)_C:= \{(x,y):\ f(x) < y < g(x)\},$$ $$(-\infty,f)_C:= \{(x,y):\ y < f(x)\}$$ and $$(f,+\infty)_C:= \{(x,y):\ y>f(x)\},$$ where $C \subseteq \mathbb R^{n-1}$ is a cell definable from $\mathcal L$ and $f,g:C \longrightarrow \mathbb R$ are continuous functions such that $f(x) < g(x)$ for $x \in C$ and the graphs of $f$ and $g$ are definable from $\mathcal L$.\medskip
	
	\noindent\textbf{Theorem (Pillay and Steinhorn, see \cite[Chapter 3]{Dries:1998nj}).} \textsl{If an $\mathcal L$-structure is o-minimal, then every set definable from $\mathcal L$ is a finite union of cells definable from $\mathcal L$.}\medskip
	
	The inductive definition of ``cell'' implies that if $C \subseteq \mathbb R^{m+n}$ is a cell definable from $\mathcal L$ and $\mu \in \mathbb R^m$, then the \textit{fiber} $C_\mu:= \{x \in \mathbb R^n:\ (\mu,x) \in C\}$ of $C$ over $\mu$ is also a cell definable from $\mathcal L$.  Thus, if $A \subseteq \mathbb R^{m+n}$ is a union of $N$ cells, where $N \in \mathbb N$, then for every $\mu \in \mathbb R^m$, the fiber $A_\mu$ is a union of at most $N$ cells.  Since the only finite cells are points, it follows that:\medskip  
	
	\noindent\textbf{Corollary (uniform finiteness principle).}
		 \textsl{If an $\mathcal L$-structure is \hbox{o-minimal}, $A \subseteq \mathbb R^{m+n}$ is definable from $\mathcal L$ and the fiber $A_\mu$ is finite for every $\mu \in \mathbb R^m$, then there exists an $N \in \mathbb N$ such that each $A_\mu$ has at most $N$ elements.}
	
\subsection*{Back to Roussarie's conjecture}
	
	One might apply the uniform finiteness principle to the set $A$ of Example 4 as follows:  let $\mathcal L_{\text{trans}}$ be the language containing $\mathcal L_{\text{or}}$ as well as the parametric transition maps associated to every limit periodic set of every $F_\mu$ in $\mathcal S_d$ as above.  Let $\mu$ be a parameter and $\Gamma$ a limit periodic set of $F_\mu$; by Example 4, the corresponding set $A$ is definable from $\mathcal L_{\text{trans}}$, and by Dulac's problem, each fiber $A_\mu$ is finite.  Therefore, the corresponding finite cyclicity conjecture conjecture follows from the uniform finiteness principle and the following:\medskip
	
	\noindent\textbf{Conjecture (o-minimality).} \textsl{The $\mathcal L_{\text{trans}}$-structure on the real numbers is \hbox{o-minimal}.}\medskip
	
	This conjecture is open, and proving o-minimality of  \hbox{$\mathcal L$-structures} is a long process.  However, a few general methods for doing so are now established and have been successfully used to obtain the following special case of the o-minimality conjecture: let $\mathcal {NRH}_d$ be the subset of all vector fields in $\mathcal S_d$ that have only \textit{non-resonant hyperbolic} singularities, as defined in the introduction of Kaiser et al. \cite{Kaiser:2009ud}.  Let $\mathcal L_{\text{nrhyp}}$ be the sublanguage of $\mathcal L_{\text{trans}}$ consisting of all parametric transition maps associated to the vector fields in $\mathcal{NRH}_d$.  Then:\medskip
	
	\noindent\textbf{Theorem \cite{Kaiser:2009ud}.}
		\textsl{The $\mathcal L_{\text{nrhyp}}$-structure on the real numbers is o-minimal; in particular, Roussarie's conjecture holds for $\mathcal {NRH}_d$.}
	\medskip
	
	The set $\mathcal {NRH}_d$ is arguably a very ``small'' subset of $\mathcal S_d$; for instance, it is not \textit{generic}, which means that even if $F_\mu \in \mathcal {NRH}_d$ for some $\mu$, there are arbitrarily close $\mu'$ such that $F_{\mu'}$ belongs to $\mathcal{S}_d$, but not to $\mathcal {NRH}_d$.  However, the larger set $\mathcal H_d$ of all vector fields in $\mathcal S_d$ that have only \textit{hyperbolic} singularities (including resonant ones) is a generic subset of $\mathcal S_d$.  One simplification for the subfamily $\mathcal{H}_d$ over the general case is that every limit periodic set is indeed a polycycle in this situation (this follows from \cite[Theorem 5 of Chapter 2]{MR3134495}, because hyperbolic singularities are always isolated).  In collaboration with my former student Zeinab Galal and my colleagues Tobias Kaiser, Jean-Philippe Rolin and Tamara Servi, I am currently working on the \hbox{o-minimality} conjecture for the corresponding sublanguage of $\mathcal L_{\text{trans}}$.  

\end{document}